\documentclass[english,11pt]{article}

\usepackage[latin1]{inputenc}
\usepackage{amsmath,amsthm,amssymb}
\usepackage{graphicx}
\usepackage{color}

\def\proof{{\it Proof.}\ }
\def\proofoft{{\it Proof of the theorem.}\ }
\def\endproof{\nolinebreak\hfill $\square$ \par\vskip3mm}

\def\eq#1{(\ref{#1})}

\def\neweq#1{\begin{equation}\label{#1}}
\def\endeq{\end{equation}}

\def\ep{\varepsilon}
\def\la{\lambda}

\def\vp{\varphi_1}
\def\RR{{\mathbb R} }

\def\di{\displaystyle}
\def\ri{\rightarrow}

\def\ii{\^\i }
\newtheorem{theorem}{Theorem}[section]
\newtheorem{lem}{Lemma}[section]
\newtheorem{corollary}{Corollary}[section]

\title{\sc Singular elliptic problems with convection term
 in anisotropic media}
\author{Louis DUPAIGNE$^a$, Marius GHERGU$^b$ and
Vicen\c tiu R\u ADULESCU$^b$\\
\small $^a$ LAMFA, Facult\'e de Math\'ematiques et d'Informatique,
Universit\'e de Picardie Jules Verne,\\
\small   33 rue Saint-Leu, 80039 Amiens, France. E-mail: {\tt
louis.dupaigne@u-picardie.fr}\\
\small $^b$ Department of Mathematics, University of Craiova,
200585 Craiova, Romania\\ \small E-mail: {\tt
marius.ghergu@math.cnrs.fr}
\qquad {\tt vicentiu.radulescu@math.cnrs.fr}}

\date{}

\begin{document}
\baselineskip16pt \maketitle
\renewcommand{\theequation}{\arabic{section}.\arabic{equation}}
\catcode`@=11 \@addtoreset{equation}{section} \catcode`@=12

\begin{abstract}
We are concerned with singular elliptic problems of the form
$-\Delta u\pm p(d(x))g(u)=\la f(x,u)+\mu |\nabla u|^a$ in $\Omega,$
where $\Omega$ is a smooth bounded domain in $\RR^N$, $d(x)={\rm
dist}(x,\partial\Omega),$
$\la>0,$ $\mu\in\RR$, $0<a\leq 2$, and  $f,k$ are  nonnegative and
nondecreasing functions. We
assume that $p(d(x))$ is a positive weight with possible singular
behavior
on the boundary of $\Omega$ and that the nonlinearity $g$ is
unbounded around the origin. Taking into account the competition
between the anisotropic potential $p(d(x))$, the convection term
$|\nabla u|^a$, and  the singular nonlinearity $g$, we establish
various existence and nonexistence results.
\end{abstract}

\section{Introduction}

Let $\Omega\subset \RR^N$ $(N\geq 2)$ be a bounded domain with
smooth boundary. We are concerned in this paper with singular
elliptic problems of the following type
$$\left\{\begin{tabular}{ll}
$-\Delta u\pm p(d(x))g(u)=\la f(x,u)+\mu |\nabla u|^a$
\quad & $\mbox{\rm in}\ \Omega,$\\
$u>0$ \quad & $\mbox{\rm in}\ \Omega,$\\
$u=0$ \quad & $\mbox{\rm on}\ \partial\Omega,$
\end{tabular} \right.\eqno(P)^\pm$$
where $d(x)=\mbox{dist}(x,\partial\Omega),$ $\la>0$, $\mu\in\RR$,
and $0<a\leq 2$.

We refer the reader to the works of Serrin \cite{serrinacta},
Choquet-Bruhat and Leray \cite{leray}, and Kazdan and Warner
\cite{kazdan}, which motivate the requirement that the nonlinearity
$|\nabla u|^a$ grows at most quadratically. We also assume that

\begin{itemize}

\item $g\in C^1(0,\infty)$ is a positive decreasing function and

\medskip
\noindent $\di(g1)\qquad\lim_{t\ri 0^+}g(t)=+\infty.$

\medskip

\item $f:\overline{\Omega}\times[0,\infty)\rightarrow[0,\infty)$ is a
H\"{o}lder continuous function which is nondecreasing with respect
to the second variable and such that $f$ is positive on
$\overline\Omega\times (0,\infty).$ Furthermore, $f$ is either
linear or $f$ is sublinear with respect to the second variable. This
last case means that $f$ fulfills the hypotheses
\medskip

\noindent $\di(f1)\qquad$ the mapping $\di (0,\infty)\ni
t\longmapsto\frac{f(x,t)}{t}\quad\mbox{is nonincreasing for
all}\;\, x\in\overline{\Omega};$
\medskip

\noindent $\di(f2)\qquad \lim_{t\ri 0^+}
\frac{f(x,t)}{t}=+\infty\quad\mbox{and}\;\;
\lim_{t\rightarrow+\infty}\frac{f(x,t)}{t}=0,\;\;\mbox{uniformly
for}\;\,x\in\overline{\Omega}.$

\item $p:(0,+\infty)\ri(0,+\infty)$ is
nonincreasing and H\"older continuous.
\end{itemize}
\medskip
Such singular boundary value problems arise in the context of
chemical heterogeneous catalysts and chemical catalyst kinetics, in
the theory of heat conduction in electrically conducting materials,
singular minimal surfaces, as well as in the study of non-Newtonian
fluids or boundary layer phenomena for viscous fluids (we refer for
more details to \cite{bn,caf,cn,coh, diaz,dmo} and the more recent
papers \cite{crg,gr,hai,her,her1,mea,sy1,sy2,zedimb}). We also point
out that, due to the meaning of the unknowns (concentrations,
populations, etc.), only the positive solutions are relevant in most
cases.

To the best of our knowledge, there does not exist a qualitative theory
for
the study
of singular boundary value problems with nonlinearities in the Kato
class
$K_N^{\rm loc}(\RR^N)$. This theory was introduced by Aizenman and
Simon in
\cite{aize} to describe wide classes of functions arising in Potential
Theory.
We refer to the recent paper \cite{rhouma} for existence and
bifurcation
results on Dirichlet boundary value problems with indefinite
nonlinearities.

The results in this paper complete the study developed in \cite{gre}
and \cite{grjmaa} since here we deal with  singular weights. One of
our purposes is to give a necessary and sufficient condition on the
weight $p$ in order to obtain a classical solution of problems
$(P)^\pm.$ By classical solution
 we understand a function $u\in C^2(\Omega)\cap C(\overline\Omega)$
that fulfills $(P)^\pm.$

Dealing with problem $(P)^+$ we show that a necessary condition in
order
to have classical solution is
\neweq{pu0}
\di \int_0^1 p(t)g(t) dt<+\infty.
\endeq
In the case where $f$ is sublinear, that is, $f$ fulfills the
hypotheses $(f1)$ and $(f2),$ condition
\eq{pu0} is also sufficient for existence of a classical solutions of
$(P)^+$ provided $\la$ and $\mu$ belong to a certain range (see Theorem
\ref{th2}).
Obviously, \eq{pu0} implies the following Keller-Osserman type
condition around the origin (see the
proof of Theorem \ref{th2})

\medskip
\noindent $\di (KO)\qquad \int^1_{0}\left(
\int^t_0\Phi(s)ds\right)^{-1/2}dt<+\infty\,,$
where $\Phi(s)=p(s)g(s),$ for all $s>0$.

\medskip

As proved by B\'enilan, Brezis and Crandall \cite{bbc}, condition
$(KO)$ is equivalent to the {\it property of compact support},
that is, for every $h\in L^1(\RR^N)$ with compact support, there
exists a unique $u\in W^{1,1}(\RR^N)$ with compact support such
that $\Delta u\in L^1(\RR^N)$ and
$$-\Delta u+\Phi(u)=h\qquad\mbox{a.e. in}\ \RR^N.$$

The results are completely different for problem $(P)^-.$ Our
results in this case generalize those established in \cite{zc}, in the
sense that in the present
paper we do not  prescribe the behavior of the singular
nonlinearity $g$ around the origin.  Also, we proved in
\cite{gre} that if $p\equiv 1$, then the
existence of a classical
solution to $(P)^-$ does not depend on the asymptotic behavior of
$g$ near the origin, whereas the exponent $a$ of the convection
term $|\nabla u|^a$ plays a crucial role. In our case, the potential
$p(d(x))$
also affects the existence of classical solutions to $(P)^-.$

Many papers have been devoted to the case $p\equiv 1$ and $\mu=0$ (see
\cite{cp, crt, gr, sy1}
and the references therein).
One of the first works in the literature dealing with singular weights
in
connection with
singular nonlinearities is due tu Taliaferro \cite{tal}. In \cite{tal}
the
following problem has been considered
\neweq{taliaferro}
\left\{\begin{tabular}{ll}
$-y''=\varphi(x)y^{-\beta}$
\quad & $\mbox{\rm in}\ (0,1),$\\
$y(0)=y(1)=0,$
\end{tabular} \right.
\endeq
where $\beta>0$ and $\varphi(x)$ is positive and continuous on $(0,1)$.
It was proved that problem \eq{taliaferro} has solutions if and only if
$\int_0^1 t(1-t)\varphi(t)dt<+\infty.$ Later, Agarwal and O'Regan
\cite[Section 2]{ar}
studied the more general problem
\begin{equation}\label{dunu}
\left\{\begin{tabular}{ll}
$H''(t)=-p(t)g(H(t))$&$\quad \mbox{ in } (0,1),$\\
$H>0$&$\quad\mbox{ in } (0,1),$\\
$H(0)=H(1)=0,$
\end{tabular} \right.
\end{equation}
where $g$ satisfies $(g1)$ and $p$ is positive and continuous on
$(0,1)$.
It is shown in \cite{ar} that if
\neweq{pu}
\di \int_0^1 t(1-t)p(t) dt<+\infty,
\endeq
then \eq{dunu} has at least one classical solution. In our
framework, $p$ is continuous at $t=1$ so condition \eq{pu} reduces
to
\neweq{pu1}
\di \int_0^1 tp(t) dt<+\infty.
\endeq
We prove that assumption \eq{pu1} is necessary in order that problem
$(P)^-$ has classical solutions. Furthermore, we argue in Section 3
that the existence of a classical solution of $(P)^-$ when $f$ is
sublinear depends on the asymptotic behavior of the gradient term
$|\nabla u|^a.$ In this sense, we prove that if $0<a<1$, then
 $(P)^-$ has at least one classical solution
for all $\mu\in\RR$. In turn, if  $1<a\leq 2$, then $(P)^-$ has no
solutions for large values of $\mu.$

Special attention is payed to the case where $a=1$. This case  was
left as an open question in \cite{gre}. We prove in Theorem
\ref{thklin} that if $\Omega$ is a ball centered at the origin, then
$(P)^-$ has at least one solution for all $\mu\in\RR$, provided
$a=1$.

The existence of a solution to $(P)^\pm$ is achieved by the sub
and super-solution method.
In particular, the
super-solution of $(P)^-$ is expressed in terms of $H.$ In the
 case of pure power nonlinearities,
a careful analysis of \eq{dunu} allows us to give boundary
estimates of the solution.

The outline of the paper is as follows. In Section 2 we give
existence and nonexistence results for problem $(P)^+.$
Section 3 concerns the problem $(P)^-$ in which
we discuss separately the case where $f$ is linear or sublinear. At the
end
of this Section we present, as an application of the obtained results,
the case where
$p(t)=t^{-\alpha}$ and $g(t)=t^{-\beta},$ and we give some
estimates for the solution at the boundary.
To make the results clearer, we assume that
$\la=1$ and $f$ is sublinear. Thus, problem $(P)^-$ becomes
\begin{equation}\label {special}
 \left\{\begin{tabular}{ll}
$-\Delta u=d(x)^{-\alpha} u^{-\beta}+f(x,u)+\mu |\nabla u|^a$ \quad &
${\rm in}\
\Omega,$\\
$u>0$ \quad & ${\rm in}\ \Omega,$\\
$u=0$ \quad & ${\rm on}\ \partial\Omega.$
\end{tabular} \right.
\end{equation}

\section{The problem $(P)^+$}

We first establish the following general nonexistence result related to
problem $(P)^+.$

\begin{theorem}\label{nonexpminus}
Assume that $\di \int_0^1 p(t)g(t)dt=+\infty.$
Let $\Phi:\overline\Omega\times[0,+\infty)\ri\RR$ be a H\"older
continuous function. Then the inequality boundary value problem
\neweq{npminus}
\left\{\begin{tabular}{ll}
$-\Delta u+p(d(x)) g(u)\leq \Phi(x,u)+C\,|\nabla u|^2$ \quad & $\mbox{\rm
in}\
\Omega,$\\
$u>0$ \quad & $\mbox{\rm in}\ \Omega,$\\
$u=0$ \quad & $\mbox{\rm on}\ \partial\Omega,$
\end{tabular} \right.
\endeq
 has no classical solutions.
\end{theorem}

As a direct consequence, we obtain :

\begin{corollary}\label{coroplus}
Assume that $\di \int_0^1 p(t)g(t)dt=+\infty$. Then problem $(P)^+$
has no classical solutions.
\end{corollary}

\proofoft We apply an idea found e.g. in \cite{zna}. It is readily
seen that it suffices to prove the Theorem only for the case $C>0$.
Arguing by contradiction, we assume that the boundary value
inequality problem \eq{npminus} has a solution $u\in C^2(\Omega)\cap
C(\overline\Omega)$. By the Gelfand transform $v=e^{Cu}-1$ we find
\neweq{gelfa}
\begin{aligned}
 \Delta v&=\di Ce^{Cu}\left(\Delta u+C|\nabla u|^2\right)\geq
Ce^{Cu} \left[p(d(x)) g(u)- \Phi(x,u)\right]\\
& =\di C(v+1)\left[p(d(x))g\left(\frac{\ln (v+1)}{C}\right)-
\Phi \left(x,\frac{\ln (v+1)}{C}\right)\right].\\
\end{aligned}
\endeq
Since $v$ is continuous on $\overline\Omega$ and $v>0$ in $\Omega$, we
deduce that
$$-\Delta v\leq \Psi(x,v)\leq C_0\quad\mbox{ in } \Omega,$$
where $\Psi(x,v)=C(v+1)\Phi \left(x,\frac{\ln (v+1)}{C}\right)$.
A straightforward argument based on the maximum principle combined with
the observation that $v=0$ on
$\partial\Omega$ shows that $v\leq c_0 d(x)$ in
$\Omega$.

For $\ep>0$ small enough, consider an open set $\Omega_\ep$ with smooth
boundary such
that $\Omega_\ep\supset\{x\in\Omega;\ {\rm dist}\,
(x,\partial\Omega)>\ep\}$. By integration in \eq{gelfa}
we find
$$-\int_{\partial\Omega_\ep}\frac{\partial
v}{\partial\nu_\ep}ds+C\int_{\Omega_\ep}(v+1)
p(d(x))g\left(\frac{\ln (v+1)}{C}\right)dx\leq
\int_{\Omega_\ep}\Psi(x,v)dx\leq C_0|\Omega|.$$
Therefore
\neweq{constp00}
\begin{aligned}
\int_{\partial\Omega_\ep}\frac{\partial
v}{\partial\nu_\ep}ds&\geq C\int_{\Omega_\ep}(v+1)
p(d(x))g\left(\frac{\ln (v+1)}{C}\right)dx-
C_0|\Omega|\\
&\geq C\int_{\Omega_\ep}
p(d(x))g\left(\frac{v}{C}\right)dx-
C_0|\Omega|.\\
\end{aligned}
\endeq
Since $v\leq c_0 d(x)$ in
$\Omega$, and $\int_{0}^1p(t)g(t)dt=+\infty$, it follows that the integral in the right-hand side of
\eq{constp00} diverges as $\ep\ri 0^+$. Hence
$$\di \lim_{\ep\ri 0^+}\int_{\partial\Omega_\ep}\frac{\partial
v}{\partial\nu_\ep}ds=+\infty.$$
But this contradicts the maximum principle (see \cite[Lemma 3.4]{gt})
because $\limsup_{t\ri 0^-}\frac{v(x_0+t\nu)}{t}<0$, for all $x_0\in\partial\Omega$.
\endproof

Before stating our existence results, we recall the following
auxiliary tool (see \cite[Lemma 2.1]{grjmaa} for a complete proof).

\begin{lem}\label{l1}
Let $\Psi:\overline{\Omega}\times(0,+\infty)\rightarrow\RR$ be a
H\"older continuous function such that the mapping
$\di(0,+\infty)\ni s\longmapsto\frac{\Psi(x,s)}{s}$ is strictly
decreasing for each $x\in\Omega.$ Assume that there exist  $v$,
$w\in C^2(\Omega)\cap C({\overline{\Omega}})$ such that

$(a)\qquad \Delta w+\Psi(x,w)\leq 0\leq \Delta v+\Psi(x,v)$ in
$\Omega;$

$(b)\qquad v,w>0$ in $\Omega$ and $v\leq w$ on $\partial\Omega;$

$(c)\qquad\Delta v\in L^1(\Omega)\;\mbox{ or }\;
\Delta w\in L^1(\Omega).$

Then $v\leq w$ in $\Omega.$
\end{lem}

Next, we prove that \eq{pu0} is sufficient for the existence of a
classical solution to $(P)^+$ provided $\mu\leq 0$ and $\la>0$ is sufficiently large. We have

\begin{theorem}\label{th2}
Assume that $\di \int^1_0 p(t)g(t)dt<+\infty.$
\begin{enumerate}
\item[{\rm (i)}] If $\mu=-1,$ then there exists $\la^*>0$ such that
$(P)^+$ has at
least one classical solution if $\la>\la^*$ and no solution exists if
$\;0<\la<\la^*.$
\item[{\rm (ii)}] If $\mu=+1$ and $0<a<1$, then there
exists $\la^*>0$ such that $(P)^+$ has at least one classical
solution for all $\la>\la^*$ and no solution exists if
$\;0<\la<\la^*.$
\end{enumerate}
\end{theorem}

\proof (i) We split the proof into several steps.\\
\noindent{\it Step 1: Existence of a solution for $\la$ large.}
By virtue of \cite[Lemma 4]{sy1} (see also (\cite[Theorem 2]{sy2}),
the problem
\neweq{UU}
\left\{\begin{tabular}{ll}
$-\Delta U=\la f(x,U)$ \quad & $\mbox{\rm in}\ \Omega,$\\
$U>0$ \quad & $\mbox{\rm in}\ \Omega,$\\
$U=0$ \quad & $\mbox{\rm on}\ \partial\Omega,$
\end{tabular} \right.
\endeq
has at least one classical solution $U_\la$, for all $\la>0$.
Using the regularity of
$f$ it
follows that $U_\la \in C^2(\overline\Omega)$ and there
exist
$c_1,c_2>0$ depending on $\la$ such that
\neweq{uuc}
c_1 d(x)\leq U_\la (x) \leq c_2 d(x)\quad\mbox{ in }\;\Omega.
\endeq

Fix $\la>0$ and observe that $U_\la$ is a super-solution of
$(P)^+.$ The main point is to find a
sub-solution $\underline u_\la$ of $(P)^+$ such that $\underline
u_\la\leq
U_\la$ in $\Omega.$
For this purpose, let
$\Phi(t)=p(t)g(t),$ $t>0$, and define
$$\di \Psi:[0,+\infty)\ri[0,+\infty),\quad
\Psi(t)=\int_0^t\frac{1}{\sqrt{2\int_0^s\Phi(\tau)d\tau}}ds.$$
Remark first that $\Psi$ is well defined, since $\Phi\in
L^1(0,1).$ Indeed, there exists $m>0$ such that
$\Phi(s)\geq m,$ for all $0<s<1.$ This  yields
$(\int_0^s\Phi(\tau)d\tau)^{-1/2}\leq (\sqrt{ms})^{-1},$
for all $0<s<1$ which implies the Keller-Osserman condition $(KO)$
around the origin:
$$\di \int^1_{0}\left(\int^t_0\Phi(s)ds\right)^{-1/2}dt<+\infty.$$

We claim that $\Psi$ is a bijective map. Indeed, $\Psi$
is increasing and
if $M:=\Phi(1),$ then
$$\di \int_0^s\Phi(\tau)d\tau\leq
\int_0^1\Phi(\tau)d\tau+M(s-1),\quad\forall s\geq 1.$$
Thus, there exists $c>0$ such that
$$\di \int_0^s\Phi(\tau)d\tau\leq Ms+c,\quad\forall s\geq 1.$$
It follows that
$$\di \Psi(t)\geq \int_1^t\frac{1}{\sqrt{2(Ms+c)}}ds\geq
\frac{1}{M}(\sqrt{2(Mt+c)}-c_1),\quad\forall\,t\geq 1.$$
This gives $\lim_{t\ri+\infty}\Psi(t)=+\infty$
and the claim follows.

Let $h:[0,+\infty)\rightarrow[0,+\infty)$ be the inverse of $\Psi.$ Then
$h$ satisfies
\begin{equation}\label{doicinci}
\left\{\begin{tabular}{ll}
$h>0$&$\quad\mbox{ in }\;(0,+\infty),$\\
$h'(t)=\sqrt{2\int_0^{h(t)}\Phi(s)ds}$&$ \quad\mbox{ in
}\;(0,+\infty),$\\
$h''(t)=\Phi(h(t))$&$\quad\mbox{ in }\;(0,+\infty),$\\
$\di h(0)=h'(0)=0.$
\end{tabular} \right.
\end{equation}
Hence $h\in C^2(0,+\infty)\cap C^1[0,+\infty).$
Let $\varphi_1>0$ be the first eigenfunction of $(-\Delta)$
in $H^1_0(\Omega)$. It is well known that there exists $C>0$ such that
\neweq{constp}
 C d(x)\leq \vp\leq\frac{1}{C}d(x)\quad\text{ for all }\,x\in\Omega.
\endeq

The key result for this part of the proof is the following.

\begin{lem}\label{lstep}
There exist two positive constants $c>0$ and $M>0$ such that
$\underline{u}_{\la}:=Mh(c\vp)$
is a sub-solution of $(P)^+$ provided $\la>0$ is large enough.
\end{lem}

\proof Since $h\in C^1[0,\infty)$ and $h(0)=0$,
we can take $c>0$ small enough such that
\neweq{himp}
h(c\vp)\leq d(x) \quad\mbox{ in }\ \Omega.
\endeq
By Hopf's maximum principle, there exist $\delta>0$ and
$\omega\subset\subset\Omega$ such that $\di |\nabla \vp|\geq\delta$
in $\Omega\setminus\omega.$ Let
\neweq{mm1}
\di M=\max\{1,2(c\delta)^{-2}\}.
\endeq
Since
$$\lim_{d(x)\ri 0^+}
\Big\{-p(d(x))g(h(c\vp))+Mc\la_1\vp
h'(c\vp)+(Mch'(c\vp)|\nabla\vp|)^a\Big\}=-\infty,$$
we can assume that
\neweq{omegaz}
-p(d(x))g(h(c\vp))+Mc\la_1\vp
h'(c\vp)+(Mch'(c\vp)|\nabla\vp|)^a<0\quad\mbox{ in }\Omega\setminus\omega.
\endeq
We are now able to show that $\underline{u}_{\la}:=Mh(c\vp)$
is a sub-solution of $(P)^+$
provided $\la>0$ is sufficiently large. Indeed, using the monotonicity
of
$g$ and \eq{himp} we have
\begin{equation}\label{doinoua}
\begin{tabular}{lll}
$\di -\Delta
\underline{u}_{\la}+p(d(x))g(\underline{u}_{\la})+|\nabla\underline
u_\la|^a=$\\
$\quad=-Mc^2p(h(c\vp))g(h(c\vp))|\nabla\vp|^2+Mc\la_1\vp
h'(c\vp)$\\
$\quad\di
\;\;\;\,+p(d(x))g(Mh(c\vp))+(Mch'(c\vp)|\nabla\vp|)^a$\\
$\quad\leq p(d(x))g(h(c\vp))(1-Mc^2|\nabla\vp|^2)+Mc\la_1\vp
h'(c\vp)+(Mch'(c\vp)|\nabla\vp|)^a.$
\end{tabular}
\end{equation}
Taking into account the definition of $M$ and \eq{omegaz}, we find
\neweq{bord}
\begin{tabular}{ll}
$\di -\Delta
\underline{u}_{\la}+p(d(x))g(\underline{u}_{\la})+(|\nabla\underline
u_\la|)^a$\\
$\qquad \leq -p(d(x))g(h(c\vp))+Mc\la_1\vp
h'(c\vp)+(Mch'(c\vp)|\nabla\vp|)^a<0\quad\mbox{ in }\Omega\setminus\omega.$
\end{tabular}
\endeq
On the other hand, from \eq{doinoua} and for all $x\in\omega$ we have
\neweq{ddoix}
\begin{aligned}
\di -\Delta
\underline{u}_{\la}+p(d(x))g(\underline{u}_{\la})+|\nabla\underline
u_\la|^a
\leq&\, p(d(x))g(h(c\vp))+Mc\la_1\vp
h'(c\vp)\\
&+(Mch'(c\vp)|\nabla\vp|)^a.\\
\end{aligned}
\endeq
Since $\vp>0$ in $\overline\omega$ and $f$ is positive on
$\overline \omega\times(0,+\infty),$ we may choose $\la>0$ such
that
\neweq{lambda2}
\di \la\min_{x\in\overline\omega} f(x,Mh(c\vp))\geq
\max_{x\in\overline\omega}\Big\{p(d(x))g(h(c\vp))+Mc\la_1\vp
h'(c\vp)+ (Mch'(c\vp)|\nabla\vp|)^a\Big\}.
\endeq
From \eq{ddoix}  and \eq{lambda2} we deduce
\neweq{intfinal}
\di -\Delta
\underline{u}_{\la}+p(d(x))g(\underline{u}_{\la})+|\nabla\underline
u_\la|^a \leq \la f(x,\underline u_\la)\quad\mbox{ in }\omega.
\endeq
Now, relations \eq{bord} and \eq{intfinal} show that $\underline
u_\la=Mh(c\vp)$ is a sub-solution of $(P)^+$ provided $\la>0$
satisfies \eq{lambda2}. This finishes the proof of our Lemma.
\endproof

Using Lemma \ref{l1}, it follows that
$\underline u_\la\leq U_\la$ in $\Omega$ and by standard elliptic
arguments (see \cite{gt}) we obtain a classical solution $u_\la$ of
$(P)^+$
such that $\underline u_\la\leq u_\la\leq U_\la$ in $\Omega.$

\medskip
\noindent {\it Step 2: Nonexistence for $\la>0$ small.}
We first remark that
$$\di \lim_{t\ri 0^+}(f(x,t)-p(d(x))g(t))=-\infty \quad
\mbox{ uniformly for }\, x\in\Omega.$$
Hence, there exists $t_0>0$ such that
\neweq{n1}
\di f(x,t)-p(d(x))g(t)<0, \quad \mbox{ for all }\,
(x,t)\in\Omega\times(0,t_0).
\endeq
On the other hand, the assumption $(f1)$ yields
\neweq{n2}
\di \frac{f(x,t)-p(d(x))g(t)}{t}\leq \frac{f(x,t)}{t}\leq
\frac{f(x,t_0)}{t_0},
\endeq
for all $(x,t)\in\Omega\times[t_0,+\infty).$
Let $m=\max_{x\in\overline\Omega}\frac{f(x,t_0)}{t_0}.$ Combining
\eq{n1} with \eq{n2} we find
\neweq{n3}
\di f(x,t)-p(d(x))g(t)<mt, \quad \mbox{ for all }\,
(x,t)\in\Omega\times(0,+\infty).
\endeq
Set $\la_0=\min\left\{1,\la_1/2m\right\}.$ We claim that problem
$(P)^+$ has no classical solution for $0<\la\leq \la_0.$ Indeed,
assume by contradiction that $u_0$ is a classical solution of
$(P)^+$ with $\la\in(0,\la_0].$ Then, according to \eq{n3}, $u_0$ is
a sub-solution of
\neweq{n4}
\left\{\begin{tabular}{ll} $\di-\Delta u=\frac{\la_1}{2}u$ \quad &
${\rm in}\
\Omega,$\\
$u>0$ \quad & ${\rm in}\ \Omega,$\\
$u=0$ \quad & ${\rm on}\ \partial\Omega.$
\end{tabular} \right.
\endeq
By Lemma \ref{l1} we have
$u_0\leq U_\la$ in $\Omega.$ Furthermore, from \eq{constp} and \eq{uuc}
we get $c u_0\leq \vp$ in $\Omega$ for some positive constant $c>0.$
Note that $cu_0$ is still
a sub-solution of \eq{n4} while $\vp$ is a super-solution of \eq{n4}.
By standard elliptic arguments, problem \eq{n4} has a solution
$u\in
C^2(\overline\Omega).$ Multiplying by $\vp$ in \eq{n4} and
integrating on $\Omega$ we have
$$\di -\int_\Omega\vp\Delta udx=\frac{\la_1}{2}\int_\Omega
u\vp dx,$$ that is,
$$\di \la_1\int_\Omega u\varphi_1 dx =-\int_\Omega u
\Delta\vp dx=\frac{\la_1}{2}\int_\Omega u\vp dx.$$
The above equality yields $\int_\Omega u\vp dx=0,$ but
this is clearly a contradiction, since $u$ and $\vp$ are both
positive on $\Omega.$ It follows that $(P)^+$ has no classical
solutions for $0<\la\leq \la_0.$

\medskip
\noindent {\it Step 3: Dependence on $\la>0$.}
Set
$$ \di A=\left\{\la>0; \mbox{ problem }(P)^+
\mbox{ has at least one classical solution}\right\}.$$
From the above
arguments we deduce that $A$ is nonempty and $\la^*:=\inf A$ is
positive. We show that if $\la\in A,$ then $(\la,+\infty)\subseteq
A.$ To this aim, let $\la_1\in A$ and $\la_2>\la_1.$ If $u_{\la_1}$
is a solution of $(P)^+$ with $\la=\la_1,$ then $u_{\la_1}$ is a
sub-solution
of $(P)^+$ with $\la=\la_2$ while $U_{\la_2}$  defined in \eq{UU}
for $\la=\la_2$  is a super-solution. Moreover, we have
$$\di \Delta U_{\la_2}+\la_2f(x,U_{\la_2})\leq 0\leq \Delta u_{\la_1}+
\la_2f(x,u_{\la_1})\quad\mbox{ in }\Omega,$$
$$\di U_{\la_2},u_{\la_1}>0\quad\mbox{ in }\Omega,$$
$$U_{\la_2}=u_{\la_1}=0\quad\mbox{ on }\partial\Omega,$$
$$\Delta U_{\la_2}\in L^1(\Omega).$$
Again by Lemma \ref{l1} we get $u_{\la_1}\leq U_{\la_2}$ in
$\Omega.$ Therefore, problem $(P)^+$ with $\la=\la_2$ has at least one
classical solution.  Since $\la\in A$ was arbitrary, we
conclude that $(\la^*,+\infty)\subset A.$
This completes the proof of (i).
\medskip

(ii) \noindent{\it Step 1: Existence of a solution for $\la$ large.}

According to Lemma \ref{lstep}, there exists $\la^*>0$ such that
$(P)^+$ has a sub-solution $\underline u_\la$ for $\la>\la^*$ and
$\mu=-1.$ Then $\underline u_\la$ is also a sub-solution in case
$\mu=+1$, provided $\la>\la^*.$ Let us construct now a
super-solution. By \cite[Lemma 4]{sy1}, for all $\la>\la^*$ there
exists $v_\la\in C^2(\overline\Omega)$ a solution of
$$
\left\{\begin{tabular}{ll}
$-\Delta v=\la f(x,v)+1$ \quad & $\mbox{\rm in}\ \Omega,$\\
$v>0$ \quad & $\mbox{\rm in}\ \Omega,$\\
$v=0$ \quad & $\mbox{\rm on}\ \partial\Omega.$
\end{tabular} \right.
$$
Since  $0<a<1$, we can choose $M=M(\la)>1$ large enough such that
$M>M^a|\nabla v_\la|^a$ in $\Omega$. Then, using $(f1)$ we obtain
$$\di -\Delta(M v_\la)=\la Mf(x,v_\la)+M\geq \la f(x,Mv_\la)+|\nabla
(Mv_\la)|^a\quad\text{ in }\Omega.$$ Hence $\overline
u_\la:=Mv_\la\in C^2(\overline\Omega)$ is a super-solution of
$(P)^+$ for all $\la>\la^*.$ On the other hand, since $\Delta
\overline u_\la+\la f(x,\overline u_\la)\leq 0\leq \Delta \underline
u_\la+\la f(x,\underline u_\la)$ in $\Omega,$ by Lemma \ref{l1} we
get $\underline u_\la\leq \overline u_\la$ and finally, problem
$(P)^+$ has at least one solution for all $\la>\la^*.$

\medskip
\noindent {\it Step 2: Nonexistence for $\la>0$ small.} We first
extend Lemma \ref{l1} in the following way :
\begin{lem}\label{l1bis}
Let $0<a<1$ and
$\Psi:\overline{\Omega}\times(0,+\infty)\rightarrow\RR$ be a H\"older
continuous function such that the mapping $\di(0,+\infty)\ni
s\longmapsto\frac{\Psi(x,s)}{s}$ is strictly decreasing for each
$x\in\Omega.$ Assume that there exist  $v$, $w\in C^2(\Omega)\cap
C({\overline{\Omega}})$ such that

$(a)\qquad \Delta w+\Psi(x,w)+|\nabla w|^a\leq 0\leq \Delta
v+\Psi(x,v)+|\nabla v|^a$ in $\Omega;$

$(b)\qquad v,w>0$ in $\Omega$ and $v<w$ on $\partial\Omega.$

Then $v\le w$ in $\Omega.$
\end{lem}
\proof Assume by contradiction that the inequality $v\le w$ does not
hold throughout $\Omega$ and let $\varphi=\frac{v}{w}$. Clearly,
$\varphi<1$ on $\partial\Omega$ and
$$
-\nabla\cdot\left[w^2\nabla\varphi\right] = -w\Delta v + v\Delta w.
$$
Let $x_0\in\Omega$ denote a point of maximum of $\varphi$. In
particular $\nabla\varphi(x_0)=0$, $-\Delta\varphi(x_0)\ge0$ and it
follows that
$$
0\le [-w\Delta v + v\Delta w](x_0).
$$
Since $w(x_0)< v(x_0)$, it follows from assumption $(a)$, the
properties of $\Psi$ and the above inequality that
$$
0< \left[|\nabla v|^a w - |\nabla w|^a v\right](x_0).
$$
Since $\nabla\varphi(x_0)=0$, we finally obtain
$$
0< \left[\left(\frac{v}{w}\right)^a w - v\right]|\nabla w|^a
(x_0)=v^a\left(w^{1-a}-v^{1-a}\right)|\nabla w|^a(x_0),
$$
contradicting $w(x_0)< v(x_0)$. \qed

\medskip
\noindent Next, we assume by contradiction that there exists a
sequence of solutions $u_n$ of $(P^+)$ associated to a parameter
$\lambda_n\ri 0^+$. A simple calculation shows that
$w(x)=A(R^2-|x|^2)$ is positive and satisfies the inequality $\Delta
w+ f(x,w)+|\nabla w|^a\leq 0$ in $\Omega$, where $A,R>0$ are large
constants. In particular, it follows from Lemma \ref{l1bis} that
$0<u_n\le w$ whenever $\lambda_n\le1$. Let $x_n\in\Omega$ denote a
maximum point of $u_n$. Then $\nabla u_n(x_n)=0$ and $-\Delta
u_n(x_n)\ge0$. Letting $d_n=d(x_n)$, $M_n=u_n(x_n)$, it follows
from $(P^+)$ that
$$
p(d_n)g(M_n)\le \lambda_n f(x_n,M_n)\le C\lambda_n,
$$
which yields a contradiction as $n\to\infty$.

The rest of the proof of (ii) follows as in the case $\mu=-1$ and
Theorem \ref{th2} is now complete.
\endproof

\section{The problem $(P)^-$}

\subsection{A nonexistence result}

We first prove :
\begin{theorem}\label{nonexpplus}
Assume that $\di \int_0^1 t p(t)dt=+\infty$.
Then the inequality boundary value problem
\neweq{nminus}
\left\{\begin{tabular}{ll}
$-\Delta u+C|\nabla u|^2\geq p(d(x)) g(u)$ \quad & $\mbox{\rm in}\
\Omega,$\\
$u>0$ \quad & $\mbox{\rm in}\ \Omega,$\\
$u=0$ \quad & $\mbox{\rm on}\ \partial\Omega,$
\end{tabular} \right.
\endeq
 has no classical solutions.
\end{theorem}
As a direct consequence, we obtain :

\begin{corollary}\label{corominus}
Assume that $\di \int_0^1 tp(t)dt=+\infty$. Then the problem $(P)^-$
has no classical solutions.
\end{corollary}

\proofoft It suffices to prove the Theorem only for $C>0$. We argue
by contradiction and assume that there exists $u\in C^2(\Omega)\cap
C(\overline\Omega)$ a solution of \eq{nminus}. Using $(g1)$, we can
find $c_1>0$ such $\underline u:=c_1\varphi_1$ verifies
$$\di -\Delta \underline u+C|\nabla \underline u|^2\geq p(d(x)) g(\underline u) \quad\mbox{ in }\,\Omega.$$
Since $g$ is decreasing, we easily obtain
\neweq{undom}
u\geq \underline u\quad\mbox{ in }\,\Omega.
\endeq

We make in \eq{nminus} the change of variable $v=1-e^{-Cu}$. Therefore
\neweq{pvminus}
\left\{\begin{tabular}{ll}
$\di -\Delta v=C(1-v)\left(C|\nabla u|^2-\Delta u\right)\geq
C(1-v)p(d(x))g\left(
-\frac{\ln (1-v)}{C}\right)$ \quad & $\mbox{\rm in}\ \Omega,$\\
$v>0$ \quad & $\mbox{\rm in}\ \Omega,$\\
$v=0$ \quad & $\mbox{\rm on}\ \partial\Omega.$
\end{tabular} \right.
\endeq
In order to avoid the singularities in \eq{pvminus} let us consider
the approximated problem
\neweq{pvminuse}
\left\{\begin{tabular}{ll}
$\di -\Delta v=C(1-v)p(d(x))g\left(\ep
-\frac{\ln (1-v)}{C}\right)$ \quad & $\mbox{\rm in}\ \Omega,$\\
$v>0$ \quad & $\mbox{\rm in}\ \Omega,$\\
$v=0$ \quad & $\mbox{\rm on}\ \partial\Omega,$
\end{tabular} \right.
\endeq
with $0<\ep<1$.
Clearly $v$ is a super-solution of \eq{pvminuse}. Furthermore, by \eq{undom} and
the fact that $\lim_{t\ri 0^+}\frac{1-e^{-Ct}}{t}=C>0$, there exists $c_2>0$ such that
$v\geq c_2\varphi_1$ in $\Omega$. On the other hand, there exists $0<c<c_2$ such that $c\varphi_1$
is a sub-solution of \eq{pvminuse} and obviously $c\varphi_1\leq v$ in $\Omega$. Then,
the problem \eq{pvminuse} has a solution $v_\ep\in C^2(\overline\Omega)$ such that
\neweq{vpshh}
\di c\varphi_1\leq v_\ep\leq v\quad\mbox{ in }\,\Omega.
\endeq
Multiplying by $\varphi_1$ in \eq{pvminuse} and integrating we find
$$\di \lambda_1\int_\Omega \varphi_1v_\ep dx=
C\int_\Omega (1-v_\ep)\varphi_1p(d(x))g\left(
\ep-\frac{\ln (1-v_\ep)}{C}\right)dx.$$
Using \eq{vpshh} we obtain
\neweq{phicontrad}
\begin{aligned}
\di M =:\lambda_1\int_\Omega \varphi_1vdx
&\di \geq
C\int_\Omega (1-v)\varphi_1p(d(x))g\left(
-\frac{\ln (1-v)}{C}\right)dx\\
&\di \geq C_1\int_{\Omega_\delta} \varphi_1p(d(x))dx,\\
\end{aligned}
\endeq
where $\Omega_\delta\supset\{x\in\Omega;\ d(x)<\delta\}$, for some $\delta>0$
sufficiently small.
Since $\varphi_1(x)$ behaves like $d(x)$ in $\Omega_\delta$ and $\int_0^1 t
p(t)dt=+\infty$, by \eq{phicontrad} we find a contradiction. Hence,
problem \eq{nonexpplus} has no classical solutions and the proof is now complete.
\endproof
\subsection{Existence results for $(P)^-$ in the sublinear case on $f$}

Our aim here is to give existence results concerning
$(P)^-$ in case where $f$ is sublinear. Nevertheless, we prove that
condition \eq{pu1} suffices to guarantee the existence of a
classical solution for $\mu$ belonging to a certain range.

In this case the existence of a solution is strongly dependent on the
exponent $a$.
To better understand this dependence, we assume $\la=1$ but the same
results hold for any $\la>0$
(note only that the bifurcation point $\mu^*$ in the following theorem
is dependent on $\la$).

\begin{theorem}\label{th4} Assume $\la=1$, $\int_0^1
tp(t)dt<+\infty$ and conditions $(f1)$, $(f2)$, $(g1)$ and $0<a\leq 2$
are fulfilled.
\begin{enumerate}
\item[{\rm (i)}] If  $0<a<1$, then problem
$(P)^-$ has at
least one solution,
for all $\mu\in\RR$;
\item[{\rm (ii)}] If $1<a\leq 2$, then there
exists
$\mu^*>0$
such that $(P)^-$ has at least one classical
solution for all $\mu<\mu^*$ and no solution exists if
$\mu>\mu^*.$
\end{enumerate}
\end{theorem}
As a direct consequence, we obtain the following corollary, which
can be compared to Theorem \ref{th2} :
\begin{corollary}
Assume $\mu=\pm 1$, $\int_0^1 tp(t)dt<+\infty$ and conditions
$(f1)$, $(f2)$, $(g1)$ and $0<a\leq 2$ are fulfilled.
\begin{enumerate}
\item[{\rm (i)}] If  $0<a<1$, then problem
$(P)^-$ has at least one solution, for all $\lambda>0$;
\item[{\rm (ii)}] If  $<1<a\leq 2$ and $\mu=-1$, then problem
$(P)^-$ has at least one solution, for all $\lambda>0$;
\item[{\rm (iii)}] If $1<a\leq 2$ and $\mu=+1$, then there
exists $\lambda^*>0$ such that $(P)^-$ has at least one classical
solution for all $\lambda>\lambda^*$ and no solution exists if
$\lambda<\lambda^*.$
\end{enumerate}
\end{corollary}

\proofoft (i) {\sc Case $\mu> 0.$} By \cite[Lemma 4]{sy1} there
exists a classical solution $\zeta$ of the problem
\neweq{zta}
\left\{\begin{tabular}{ll}
$-\Delta \zeta=f(x,\zeta)$ \quad & $\mbox{\rm in}\ \Omega,$\\
$\zeta>0$ \quad & $\mbox{\rm in}\ \Omega,$\\
$\zeta=0$ \quad & $\mbox{\rm on}\ \partial\Omega.$
\end{tabular} \right.
\endeq
Using the regularity of $f$ we have $\zeta\in
C^2(\overline\Omega)$.
Then, $\zeta$ is a sub-solution of $(P)^-$ provided
$\mu> 0.$ We focus now on finding a super-solution $\overline u_\mu$
of $(P)^-$ such that $\zeta\leq \overline u_\mu$ in $\Omega.$

Let $H$ be the solution of \eq{dunu}.
Since $H$ is concave, there exists $H'(0+)\in (0,+\infty].$
Taking $0<b<1 $ small enough, we can assume that $H'>0$ in $(0,b],$
so $H$ is increasing on $[0,b].$ Multiplying by $H'$ in \eq{dunu}
and integrating on $[t,b],$  we find
\neweq{hprim}
\di (H')^2(t)-(H')^2(b)=2\int_t^b p(s)g(H(s))H'(s)ds\leq
2p(t)\int_{H(t)}^{H(b)} g(\tau)d\tau.
\endeq
Using the monotonicity of $g$ it follows that

\neweq{hprim2}
\di (H')^2(t)\leq 2H(b)p(t)g(H(t))+(H')^2(b), \quad\mbox { for all
}\; 0<t\leq b.
\endeq
Hence, there exist $C_1,C_2>0$ such that
\neweq{ddoi}
\di (H')(t)\leq C_1 p(t)g(H(t)),\quad\mbox{ for all }0<t\leq b
\endeq
and
\neweq{ddoibis}
\di (H')^2(t)\leq  C_2 p(t)g(H(t)),\quad\mbox{ for all }0<t\leq b.
\endeq

Now we can proceed to construct a super-solution for $(P)^-.$
First, we fix $c>0$ such that
\neweq{const1}
c\varphi_1\leq \min\{b,d(x)\}\quad\mbox{ in }\;\Omega.
\endeq
By Hopf's maximum principle, there exist
$\omega\subset\subset\Omega$ and $\delta>0$ such that
\neweq{ext}
|\nabla\varphi_1|>\delta\quad\mbox{ in }\;\Omega\setminus\omega.
\endeq
Moreover, since
$$\lim_{d(x)\ri 0^+}\left\{c^2p(c\varphi_1)g(H(c\varphi_1))
|\nabla\vp|^2-3f(x,H(c\vp))\right\}=+\infty,$$ we can assume that
\neweq{const2}
c^2p(c\varphi_1)g(H(c\varphi_1))|\nabla\vp|^2\geq
3f(x,H(c\vp))\quad\mbox{ in }\;\Omega\setminus\omega.
\endeq
Let $M>1$ be such that
\neweq{const3}
Mc^2\delta^2>3.
\endeq
Since $H'(0+)>0$ and $0<a<1$, we can choose $M>1$ such that
$$\di M\frac{(c\delta)^2}{C_1} H'(c\vp)
\geq 3\mu (McH'(c\vp)|\nabla \vp|)^a\quad\text{ in
}\,\Omega\setminus\omega,$$
where $C_1$ is the constant appearing in \eq{ddoi}. By \eq{ddoi}, \eq{ext} and
\eq{const3} we derive
\neweq{M0}
\di  Mc^2p(c\vp)g(H(c\vp))|\nabla \vp|^2\geq 3\mu
(McH'(c\vp)|\nabla \vp|)^a\quad\text{ in }\;\Omega\setminus\omega.
\endeq
Since $g$ is decreasing and $H'(c\vp)>0$ in $\overline\omega,$
there exists $M>0$ such that
\neweq{M3}
Mc\la_1\vp H'(c\vp)\geq 3p(d(x))g(H(c\vp))\quad\text{ in }\;\omega.
\endeq
In the same manner, using $(f2)$ and the fact that $\vp>0$
in $\overline\omega,$ we can choose $M>1$ large enough such that
\neweq{M1}
Mc\la_1 \vp H'(c\vp)\geq 3 \mu (MH'(c\vp)|\nabla\vp|)^a \quad\text{ in
}\;\omega,
\endeq
and
\neweq{M2}
Mc\la_1\vp H'(c\vp)\geq 3f(x,MH(c\vp))\quad\text{ in
}\;\omega.
\endeq
For $M$ satisfying \eq{const3}-\eq{M2}, we prove that
\neweq{supsol}
\overline u_\mu(x):=MH(c\varphi_1(x)), \quad \mbox{ for all
}\;x\in\Omega,
\endeq
is a super-solution of $(P)^-.$ We have
\neweq{calculp}
\di -\Delta \overline
u_\mu=Mc^2p(c\vp)g(H(c\varphi_1))|\nabla\varphi_1|^2+
Mc\la_1\varphi_1H'(c\varphi_1) \quad\mbox{ in }\;\Omega.
\endeq
We first show that
\neweq{dprima}
Mc^2p(c\vp)g(H(c\varphi_1))|\nabla\varphi_1|^2\geq
p(d(x))g(\overline u_\mu)+f(x,\overline u_\mu)+\mu |\nabla \overline
u_\mu|^a
\quad\mbox{ in }\;\Omega\setminus\omega.
\endeq
Indeed, by \eq{const1}, \eq{ext} and \eq{const3} we get
\neweq{dprimau}
\begin{aligned}
\di \frac{M}{3}c^2p(c\vp)g(H(c\varphi_1))|\nabla\varphi_1|^2&
\geq p(d(x))g(H(c\varphi_1))\\
&\geq p(d(x))g(M
H(c\varphi_1))\\
&=p(d(x))g(\overline u_\mu) \quad\mbox{ in
}\;\Omega\setminus\omega.\\
\end{aligned}
\endeq
The assumption $(f1)$ and \eq{const2} produce
\neweq{dprimat}
\begin{aligned}
\frac{M}{3}c^2p(c\vp)g(H(c\varphi_1))|\nabla\varphi_1|^2& \geq M
f(x,H(c\varphi_1))\\
&\geq f(x,MH(c\varphi_1))\\
&= f(x,\overline u_\mu)
\quad\mbox{ in }\;\Omega\setminus\omega.\\
\end{aligned}
\endeq
From  \eq{M0} we obtain
\neweq{dprimad}
\begin{aligned}
\di \frac{M}{3}c^2p(c\vp)g(H(c\varphi_1))|\nabla\varphi_1|^2
&\di\geq \mu (Mc H'(c\varphi_1)|\nabla\varphi_1|)^a\\
&=\di \mu |\nabla \overline u_\mu|^a\quad\mbox{ in
}\;\Omega\setminus\omega.
\end{aligned}
\endeq
Now, relation \eq{dprima} follows by \eq{dprimau}, \eq{dprimat} and
\eq{dprimad}.

\smallskip
Next we prove that
\neweq{ddoua}
Mc\la_1\varphi_1H'(c\varphi_1)\geq p(d(x))g(\overline
u_\mu)+f(x,\overline
u_\mu)+\mu |\nabla \overline u_\mu|^a\quad\mbox{ in }\; \omega.
\endeq
From \eq{M3} and \eq{M1} we get
\neweq{ddouau}
\begin{aligned}
\di \frac{M}{3}c\la_1\varphi_1H'(c\varphi_1)&\geq
p(d(x))g(H(c\varphi_1))\\
&\geq p(d(x))g(MH(c\varphi_1))\\
&=p(d(x))g(\overline u_\mu)\quad\mbox{ in }\; \omega\\
\end{aligned}
\endeq
and
\neweq{ddouad}
\begin{aligned}
\di \frac{M}{3}c\la_1 \varphi_1 H'(c\varphi_1)\geq&
\di \mu (Mc H'(c\varphi_1)|\nabla\varphi_1|)^a\\
=&\di \mu |\nabla
\overline u_\mu|^a \quad\mbox{ in }\; \omega.
\end{aligned}
\endeq
Finally, from \eq{M2} we derive
\neweq{ddouat}
\di \frac{M}{3}c\la_1\varphi_1H'(c\varphi_1)\geq
f(x,MH(c\varphi_1))= f(x,\overline u_\mu) \quad\mbox{ in }\;
\omega.
\endeq
Clearly, relation  \eq{ddoua} follows from
\eq{ddouau}, \eq{ddouad} and \eq{ddouat}.

Combining \eq{calculp} with \eq{dprima} and \eq{ddoua} we conclude
that $\overline u_\mu$ is a super-solution of $(P)^-.$ Thus, by
Lemma \ref{l1} we obtain $\zeta\leq \overline u_\mu$ in
$\Omega$ and by sub and super-solution method it follows that $(P)^-$
has at least one classical solution for all $\mu> 0.$

\smallskip
 {\sc Case $\mu\leq 0.$} We fix $\nu>0$ and let $u_\nu\in
C^2(\Omega)\cap C(\overline\Omega)$ be a solution of $(P)^-$ for
$\mu=\nu.$ Then $u_\nu$ is a super-solution of $(P)^-$ for all
$\mu\leq 0.$ Set
$$\di
m:=\inf_{(x,t)\in\overline\Omega\times(0,+\infty)}\Big(p(d(x))g(t)+f(x,t)\Big).$$
Since $\lim_{t\ri 0^+}g(t)=+\infty$ and the mapping $\di
(0,+\infty)\ni t\longmapsto \min_{x\in\overline\Omega}f(x,t)$ is
positive and nondecreasing, we deduce that $m$ is a positive real
number. Consider the problem
\neweq{redus}
\left\{\begin{tabular}{ll}
$-\Delta v=m+\mu |\nabla v|^a$ \quad & $\mbox{\rm in}\ \Omega,$\\
$v=0$ \quad & $\mbox{\rm on}\ \partial\Omega.$
\end{tabular} \right.
\endeq
Clearly zero is a sub-solution of \eq{redus}. Since $\mu\leq 0$, the
solution $w$ of the problem
$$
\left\{\begin{tabular}{ll}
$-\Delta w=m$ \quad & $\mbox{\rm in}\ \Omega,$\\
$w=0$ \quad & $\mbox{\rm on}\ \partial\Omega,$
\end{tabular} \right.
$$
is a super-solution of \eq{redus}.
Hence, \eq{redus} has at least one solution $v\in
C^2(\Omega)\cap C(\overline\Omega).$
We claim that $v>0$ in
$\Omega.$ Indeed, if not, we deduce that $\min_{x\in\overline\Omega}v$
is
achieved at some point $x_0\in\Omega.$ Then $\nabla v(x_0)=0$ and
$$\di -\Delta v(x_0)=m+\mu |\nabla v(x_0)|^a=m>0,\quad\mbox{
contradiction.}$$ Therefore, $v>0$ in $\Omega.$ It is easy to see
that $v$ is sub-solution of $(P)^-$ and $-\Delta v\leq
m\leq -\Delta u_\nu$ in $\Omega,$ which gives $v\leq u_\nu$ in
$\Omega.$ Again by the sub and super-solution method we conclude that
$(P)^-$ has at least one classical solution $u_\mu\in
C^2(\Omega)\cap C(\overline\Omega).$
\medskip

(ii)  The proof follows the same steps as above. The only
difference is that \eq{M0} and \eq{M1} are no more valid for any
$\mu>0.$
The main difficulty when dealing with estimates like \eq{M0} is that
$H'(c\vp)$ may
blow-up at the boundary.
However, combining the assumption $1<a\leq 2$ with \eq{ddoibis}, we can
choose $\mu>0$
small enough such that \eq{M0} and \eq{M1} hold. This implies that the
problem $(P)^-$
has a classical solution provided $\mu>0$ is sufficiently small.

Set
$$ \di A=\left\{\mu> 0; \mbox{ problem } (P)^-\mbox{ has at
least one classical solution}\right\}.$$
From the above arguments, $A$
is nonempty. Let $\mu^*=\sup A.$ We first claim that if $\mu\in
A,$ then $(0,\mu)\subseteq A.$ To this aim, let $\mu_1\in A$ and
$0< \mu_2<\mu_1.$ If $u_{\mu_1}$ is a solution of $(P)^-$ with
$\mu=\mu_1,$ then $u_{\mu_1}$ is a super-solution of $(P)^-$ with
$\mu=\mu_2$, while $\zeta$ defined in \eq{zta} is a sub-solution.
Using Lemma \ref{l1} once more, we get $\zeta\leq u_{\mu_1}$ in
$\Omega$ so $(P)^-$ has at least one classical solution for
$\mu=\mu_2.$ This proves the claim. Since $\mu_1\in A$ was
arbitrary, we conclude that $(0,\mu^*)\subset A.$

Next, we prove that
$\mu^*<+\infty.$ To this aim, we use the following result which is a
consequence of Theorem 2.1 in
\cite{ap}.

\begin{lem}\label{l2} Assume that $a>1$. Then there
exists
a positive number $\bar \sigma$ such that the problem
\begin{equation}\label {alaa}
 \left\{\begin{tabular}{ll}
$-\Delta v\geq |\nabla v|^a+\sigma$ \quad & $\mbox{\rm in}\ \Omega,$\\
$v=0$ \quad & $\mbox{\rm on}\ \partial\Omega,$
\end{tabular} \right.
\end{equation}
has no solutions for $\sigma>\bar\sigma.$
\end{lem}

Consider $\mu\in A$ and let $u_\mu$ be a classical solution of
$(P)^-.$ Set $v=\mu^{1/(a-1)}u_\mu.$ Using our assumption $1<a\leq 2$,
we deduce that $v$ fulfills
\neweq{doptispe}
\left\{\begin{tabular}{ll} $\di -\Delta v\geq |\nabla
v|^a+m\mu^{1/(a-1)}$ \quad & $\mbox{\rm in }\Omega,$\\
$v=0$ \quad & $\mbox{\rm on }\, \partial\Omega.$
\end{tabular} \right.
\endeq
According to Lemma \ref{l2}, we obtain $m\mu^{1/(a-1)}\leq
\bar\sigma,$ that is, $\di
\mu\leq \left(\frac{\bar\sigma}{m}\right)^{a-1}.$ This means that
$\di \mu^*\leq \left(\frac{\bar\sigma}{m}\right)^{a-1},$ hence
$\mu^*$ is finite. The existence of a solution in the case
$\mu\leq 0$ can be achieved in the same manner as above.

This finishes the proof of Theorem \ref{th4}.
\endproof

In what follows we discuss the case $a=1$. Note
that the method used
in Theorem \ref{th4} does not apply here for large values of
$\mu.$

Assume that $\Omega=B_R(0)$ for some $R>0,$ where
$B_R(0)=\{x\in\RR^N;\, |x|<R\}.$ In this case
and with $\la=1$,
problem $(P)^-$ becomes
\neweq{klin}
\left\{\begin{tabular}{ll}
$-\Delta u=p(R-|x|)g(u)+f(x,u)+\mu |\nabla u|$
\quad & $|x|<R,$\\
$u>0$ \quad & $|x|<R,$\\
$u=0$ \quad & $|x|=R.$
\end{tabular} \right.
\endeq

\begin{theorem}\label{thklin} Assume that $\int_0^1 tp(t)dt<+\infty.$
Then the problem \eq{klin} has at least one solution for all
$\mu\in\RR.$
\end{theorem}

\proof The case $\mu\leq 0$ is the same as in the proof of Theorem
\ref{th4} (i).
In what follows, we assume that $\mu>0$.
Using Theorem \ref{th4} (i) it is easy to see that there exists
$\underline u\in C^2(\Omega)\cap C(\overline \Omega)$ such that
$$\left\{\begin{tabular}{ll}
$-\Delta \underline u=p(R-|x|)g(\underline u)$
\quad & $|x|<R,$\\
$\underline u>0$ \quad & $|x|<R,$\\
$\underline u=0$ \quad & $|x|=R.$
\end{tabular} \right.$$
It is obvious that $\underline u$ is a sub-solution of
\eq{klin} for all $\mu>0$.
In order to provide a super-solution of \eq{klin} we consider the
problem
\neweq{klin1}
\left\{\begin{tabular}{ll}
$-\Delta u=p(R-|x|)g(u)+1+\mu |\nabla u|$
\quad & $|x|<R,$\\
$u>0$ \quad & $|x|<R,$\\
$u=0$ \quad & $|x|=R.$
\end{tabular} \right.
\endeq
We need the following auxiliary result.

\begin{lem}\label{lklin}
Problem \eq{klin1} has at least one solution.
\end{lem}

\proof
We are looking for radially symmetric solution $u$ of \eq{klin1}, that
is, $u=u(r),$ $0\leq r=|x|\leq R.$
In this case, problem \eq{klin1} becomes
\neweq{klin2}
\left\{\begin{tabular}{ll}
$\di -u''-\frac{N-1}{r}u'(r)=p(R-r)g(u(r))+1+\mu |u'(r)|$
\quad & $0\leq r<R,$\\
$u>0$ \quad & $0\leq r<R,$\\
$u(R)=0.$ \quad &
\end{tabular} \right.
\endeq
This implies
$$\di -(r^{N-1}u'(r))'\geq 0\quad\text{ for all }\;0\leq r<R,$$
which yields $u'(r)\leq 0$ for all $0\leq r<R.$ Then \eq{klin2} gives
$$\di -\left(u''+\frac{N-1}{r}u'(r)+\mu
u'(r)\right)=p(R-r)g(u(r))+1,\quad 0\leq r<R.$$
We obtain
\neweq{kl}
-(e^{\mu r}r^{N-1}u'(r))'=e^{\mu r}r^{N-1}\psi(r,u(r)),\quad 0\leq r<R,
\endeq
where
$$\di \psi(r,t)=p(R-r)g(t)+1,\quad (r,t)\in[0,R)\times(0,+\infty).$$
From \eq{kl} we get
\neweq{kl1}
u(r)=u(0)-\int^r_0 e^{-\mu t}t^{-N+1}\int^t_0 e^{\mu
s}s^{N-1}\psi(s,u(s))ds dt, \quad 0\leq r<R.
\endeq

On the other hand, in view of Theorem \ref{th4} and using the fact
that $g$ is decreasing, there exists a unique solution
$w\in C^2(B_R(0))\cap C(\overline B_R(0))$ of the problem
\neweq{klsub}
\left\{\begin{tabular}{ll}
$-\Delta w=p(R-|x|)g(w)+1$
\quad & $|x|<R,$\\
$w>0$ \quad & $|x|<R,$\\
$w=0$ \quad & $|x|=R.$
\end{tabular} \right.
\endeq
Clearly, $w$ is a sub-solution of \eq{klin1}. Due to the uniqueness and
to the symmetry of
the domain, $w$ is radially symmetric, so, $w=w(r),$ $0\leq r=|x|\leq
R.$ As above we get
\neweq{klsub1}
\di
w(r)=w(0)-\int^r_0 t^{-N+1}\int^t_0 s^{N-1}\psi(s,w(s))ds dt, \quad
0\leq r<R.
\endeq

We claim that there exists a solution $v\in C^2[0,R)\cap C[0,R]$ of
\eq{kl1} such that $v>0$ in $[0,R).$\\
Let $A=w(0)$ and define the sequence $(v_k)_{k\geq 1}$ by
\neweq{klinn}
\left\{\begin{tabular}{ll}
$\di v_k(r)=A-\int^r_0 e^{-\mu t}t^{-N+1}\int^t_0 e^{\mu
s}s^{N-1}\psi(s,v_{k-1}(s))ds dt,$ \quad& $ 0\leq r<R,\,k\geq 1,$\\
$v_0=w.$
\end{tabular}\right.\endeq
Note that $v_k$ is decreasing in $[0,R)$ for all $k\geq 0.$ From
\eq{klsub1} and \eq{klinn}
we easily check that $v_1\geq v_0$ and by induction we deduce $v_k\geq
v_{k-1}$ for all $k\geq 1.$
Hence
$$w=v_0\leq v_1\leq...\leq v_k\leq ...\leq A\quad\text{ in }\;B_R(0).$$
Thus, there exists $v(r):=\lim_{k\ri\infty} v_k(r),$ for all $0\leq
r<R$ and $v>0$ in $[0,R).$
We can now pass to the limit in \eq{klinn} in order to get that $v$ is
a solution of \eq{kl1}.
By classical regularity results we also obtain $v\in C^2[0,R)\cap
C[0,R].$ This proves the claim.

We have obtained a super-solution $v$ of \eq{klin1} such that
$v\geq w$ in $B_R(0)$. So, the problem \eq{klin1} has at least one
solution
and the proof of our Lemma is now complete.
\qed

\medskip
Let $u$ be a solution of the problem \eq{klin1}.
For $M>1$ we have
\neweq{estk}
\begin{aligned}
-\Delta (Mu)&= Mp(R-|x|)g(u)+M+\mu |\nabla(Mu)|\\
&\geq p(R-|x|)g(Mu)+M+\mu |\nabla(Mu)|.\\
\end{aligned}
\endeq
Since $f$ is sublinear, we can choose $M>1$ such that
$$M\geq f(x,M|u|_\infty)\quad\text{ in }\;B_R(0).$$
Then $\overline u_\mu:=Mu$ satisfies
$$\di -\Delta \overline u_\mu\geq p(R-|x|)g(\overline
u_\mu)+f(x,\overline u_\mu)+\mu
|\nabla\overline u_\mu|\quad\text{ in }\;B_R(0).$$
It follows that $\overline u_\mu$ is a super-solution of \eq{klin}.
Since $g$
is decreasing we easily deduce $\underline u\leq \overline u_\mu$ in
$B_R(0)$ so, problem
$(P)^-$ has at least one solution.

The proof of Theorem \ref{thklin} is now complete.
\qed

\subsection{Existence results for $(P)^-$ in the linear case on $f$}

In this section we turn to the study of problem $(P)^-$ when $f$ is
linear. More precisely, we consider the problem
\neweq{linear}
\left\{\begin{tabular}{ll}
$-\Delta u=p(d(x))g(u)+\la u+\mu |\nabla u|^a$
\quad & $\mbox{\rm in}\ \Omega,$\\
$u>0$ \quad & $\mbox{\rm in}\ \Omega,$\\
$u=0$ \quad & $\mbox{\rm on}\ \partial\Omega,$
\end{tabular} \right.
\endeq
where $\la> 0$ and $p,g$ are as in the previous sections. We assume in
what follows
that $0<a< 1$.
\smallskip

Note that the existence results established in \cite[Lemma 4]{sy1} or
\cite{sy2} do not apply here
since the mapping
$$\di \Psi(x,t)=p(d(x))g(t)+\la t,\qquad (x,t)\in \Omega\times
(0,+\infty),$$
is not defined on $\partial\Omega\times (0,+\infty)$.

\begin{theorem}\label{thlin} Assume that $\int_0^1 tp(t)dt<+\infty$ and
conditions $(g1)$, $0<a< 1$ are fulfilled.
Then for $\mu\geq 0$ the problem \eq{linear} has solutions if and only
if $\la<\la_1.$
\end{theorem}

 \proof Fix $\la\in(0,\la_1)$ and $\mu\geq 0$. By Theorem \ref{th4} (i)
there exists
$u\in C^2(\Omega)\cap C(\overline\Omega)$ a solution of the
problem
$$\left\{\begin{tabular}{ll}
$-\Delta u=p(d(x))g(u)+\mu |\nabla u|^a$
\quad & $\mbox{\rm in}\ \Omega,$\\
$u>0$ \quad & $\mbox{\rm in}\ \Omega,$\\
$u=0$ \quad & $\mbox{\rm on}\ \partial\Omega.$
\end{tabular} \right.$$
Obviously, $\underline u_{\la\mu}:=u$ is a sub-solution of \eq{linear}.
Since $\la<\la_1,$  there exists
$v\in C^2(\overline\Omega)$ such that
$$\left\{\begin{tabular}{ll}
$-\Delta v=\la v+2$
\quad & $\mbox{\rm in}\ \Omega,$\\
$v>0$ \quad & $\mbox{\rm in}\ \Omega,$\\
$v=0$ \quad & $\mbox{\rm on}\ \partial\Omega.$
\end{tabular} \right.$$
Since $0<a<1$, we can choose $M>0$ large enough such that
$$M>\la |u|_{\infty}\text{ and }\; M>\mu (M|\nabla v|)^a\text{ in
}\,\Omega.$$
Then $w:=Mv$ satisfies
$$-\Delta w\geq \la (u+w)+\mu |\nabla w|^a\quad\text{ in }\;\Omega.$$

We claim that $\overline u_{\la\mu}:=u+w$ is a super-solution of
\eq{linear}.
Indeed, we have
\neweq{ssline}
-\Delta \overline u_{\la\mu}\geq p(d(x))g(u)+\la
\overline u_{\la\mu}+\mu |\nabla u|^a+\mu |\nabla w|^a\quad\text{ in
}\;\Omega.
\endeq
Using the assumption $0<a<1$ one can easily deduce
$$\di t_1^a+t_2^a\geq (t_1+t_2)^a,\quad\text{ for all }\;t_1,t_2\geq
0.$$
Hence
\neweq{ssline2}
\di |\nabla u|^a+|\nabla w|^a
\geq (|\nabla u|+|\nabla w|)^a\geq |\nabla (u+ w)|^a\quad \text{ in
}\;\Omega.
\endeq
Combining \eq{ssline} with \eq{ssline2} we obtain
$$\di -\Delta \overline u_{\la\mu}\geq p(d(x))g(\overline
u_{\la\mu})+\la
\overline u_{\la\mu}+\mu |\nabla \overline u_{\la\mu}|^a\quad\text{ in
}\;\Omega.$$
Hence, $(\underline u_{\la\mu}, \overline u_{\la\mu})$ is an ordered
pair of sub
and super-solution of \eq{linear}, so there exists a classical solution
$u_{\la\mu}$ of \eq{linear},
provided $\mu\geq 0$ and $0<\la<\la_1.$ Assume by contradiction that
there exist
$\la\geq \la_1$ and $\mu\geq 0$ such that the problem \eq{linear} has a
classical
solution $u_{\la\mu}.$
If $m=\min_{x\in\overline\Omega} p(d(x))g(u_{\la\mu})>0$ it follows
that
$u_{\la\mu}$ is a super-solution of
\neweq{linear1}
\left\{\begin{tabular}{ll}
$-\Delta u=\la u+m$
\quad & $\mbox{\rm in}\ \Omega,$\\
$u=0$ \quad & $\mbox{\rm on}\ \partial\Omega.$
\end{tabular} \right.
\endeq
Clearly, zero is a sub-solution of \eq{linear1}, so there exists a
classical solution $u$ of \eq{linear1} such that $u\leq u_{\la\mu}$ in
$\Omega$.
By maximum principle and elliptic regularity we get $u>0$ in $\Omega$
and $u\in C^2(\overline\Omega).$ To raise a
contradiction, we proceed as in the proof of Theorem \ref{th2} (ii).

Multiplying by $\vp$ in \eq{linear1} and then integrating over $\Omega$
we find
$$\di -\int_\Omega \vp\Delta u=\la\int_\Omega u\vp+m\int_\Omega \vp.$$
This implies $\la_1\int_\Omega u\vp=\la\int_\Omega u\vp+m\int_\Omega
\vp,$ which is a contradiction, since $\la\geq \la_1$ and $m>0.$
The proof of Theorem \ref{thlin} is now complete.
\endproof

\subsection{An application}

We show here how the results in this section applies to the problem
\eq{special}.
Recall that if $\int_0^1 tp(t)dt<+\infty$ and $\mu$ belongs to a
certain range, then
Theorem \ref{th4} asserts that \eq{special} has at least one
classical solution $u_\mu$
satisfying $u_\mu\leq M H(c\vp)$ in $\Omega,$ for some $M,c>0.$ Here
$H$ is the solution of
\begin{equation}\label{dp}
\left\{\begin{tabular}{ll}
$H''(t)=-t^{-\alpha} H^{-\beta}(t),\quad \mbox{ for all } 0<t\leq
b<1,$\\
$H,H'>0\quad\mbox{\rm in}\  (0,b],$\\
$H(0)=0.$
\end{tabular} \right.
\end{equation}
With the same idea as in the proof of Theorem \ref{th4}, we can show
that there
exists $m>0$ small enough such that
$v:=mH(c\vp)$ satisfies
\neweq{inqsubb}
\di -\Delta v\leq d(x)^{-\alpha}v^{-\beta}\quad\text{ in }\;\Omega.
\endeq
Indeed, we have
$$\di -\Delta v=m[c^{2-\alpha}|\nabla \vp|^2\vp^{-\alpha}
H^{-\beta}(c\vp)+\la_1 c\vp H'(c\vp)]\quad\text{ in }\;\Omega.$$
Using \eq{constp} and \eq{ddoi}, there exist two positive constants
$c_1,c_2>0$ such that
$$\di -\Delta v\leq m[c_1|\nabla \vp|^2+c_2\vp]d(x)^{-\alpha}
H^{-\beta}(c\vp)\quad\text{ in }\;\Omega.$$
Clearly \eq{inqsubb} holds if we choose $m>0$ small enough such that
$m[c_1|\nabla \vp|^2+c_2\vp]<1$ in $\Omega.$ Moreover, $v$ is a
sub-solution of \eq{special} for all $\mu>0$
and one can easily see that
$v\leq u_\mu$ in
$\Omega.$ Hence
\neweq{asimp1}
\di mH(c\vp)\leq  u_\mu \leq M H(c\vp)\quad\mbox{ in }\;\Omega.
\endeq
Now, a careful analysis of \eq{dp} together with \eq{asimp1}
is used in order to obtain
boundary estimates for the solution of \eq{special}. Our estimates
complete the results in \cite[Theorem 2.1]{gl}
since here the potential $p(d(x))$ blows-up at the boundary.

\begin{theorem}\label{th5}
The following properties hold true.
\begin{enumerate}
\item[{\rm (i)}] If $\alpha\geq 2$, then the problem \eq{special} has
no classical
solutions.
\item[{\rm (ii)}]  If $\alpha<2$, then there exists $\mu^*\in
(0,+\infty]$
(with $\mu^*=+\infty$ if $0<a<1$) such that problem \eq{special} has at
least
one classical solution $u_\mu,$
for all $-\infty<\mu<\mu^*.$
Moreover, for all $0<\mu <\mu^*$, there exist $0<\delta<1$ and $C_1,
C_2>0$ such
that $u_\mu$ satisfies

\begin{enumerate}
\item[{\rm (ii1)}] If  $\alpha+\beta>1,$ then
\neweq{sp1}
C_1 d(x)^{\frac{2-\alpha}{1+\beta}}\leq u_\mu(x)\leq
C_2d(x)^{\frac{2-\alpha}{1+\beta}}, \quad\mbox{ for all
}\;x\in\Omega;
\endeq
\item[{\rm (ii2)}] If  $\alpha+\beta=1,$ then
\neweq{sp2}
C_1 d(x)(-\ln d(x))^{\frac{1}{2-\alpha}}\leq u_\mu(x)\leq C_2
d(x)(-\ln d(x))^{\frac{1}{2-\alpha}},
\endeq
for all $x\in\Omega$ with $d(x)<\delta;$
\item[{\rm (ii3)}] If  $\alpha+\beta<1,$ then
\neweq{sp3}
C_1d(x)\leq u_\mu(x)\leq C_2d(x), \quad\mbox{ for all }\;x\in\Omega.
\endeq
 \end{enumerate}
\end{enumerate}
\end{theorem}

\proof The existence and nonexistence of a solution to \eq{special}
follows directly from Theorems \ref{nonexpplus} and \ref{th4}.
We next prove the boundary estimates \eq{sp1}-\eq{sp3}.
\medskip

(ii1) Remark that
$$\di
H(t)=\left(\frac{(1+\beta)^2}{(2-\alpha)(\alpha+\beta-1)}\right)^{1/(1+\beta)}
t^{\frac{2-\alpha}{1+\beta}},\qquad t>0,\;$$
is a solution of
\eq{dp} provided $\alpha+\beta>1.$ The conclusion in this case
follows now from \eq{asimp1}.
\medskip

(ii2) Note that in this case problem \eq{dp} becomes
\begin{equation}\label{dunus}
\left\{\begin{tabular}{ll}
$H''(t)=-t^{-\alpha}H^{\alpha-1}(t),\quad \mbox{ for all } 0<t\leq
b<1,$\\
$H(0)=0,$\\
$H>0\quad\mbox{\rm in}\  (0,b].$
\end{tabular} \right.
\end{equation}

Since $H$ is concave, it follows that
\neweq{HH}
H(t)>tH'(t),\quad\mbox{for all }\;0<t\leq b.
\endeq
Relations \eq{dunus} and \eq{HH} yield
$$\di -H''(t)<t^{-1}(H'(t))^{\alpha-1},\quad\mbox{for all }\;0<t\leq
b.$$
Hence
\neweq{HH1}
-H''(t)(H'(t))^{1-\alpha}\leq  \frac{1}{t},\quad\mbox{for all
}\;0<t\leq b.
\endeq
Integrating in \eq{HH1} over $[t,b]$ we get
$$\di (H')^{2-\alpha}(t)-(H')^{2-\alpha}(b)\leq
(2-\alpha)(\ln b-\ln
t), \quad\mbox{for all }\;0<t\leq b.$$ Hence, there exist $c_1>0$ and
$\delta_1\in(0,b)$ such that
\neweq{HH2}
H'(t)\leq c_1(-\ln t)^{\frac{1}{2-\alpha}},\quad\mbox{for all
}\;0<t\leq \delta_1.
\endeq
Fix $t\in(0,\delta_1].$ Integrating over $[\ep,t]$, $0<\ep<t,$ in
\eq{HH2} we have
\neweq{HH3}
\di H(t)-H(\ep)\leq c_1 t(-\ln t)^{\frac{1}{2-\alpha}}+
\frac{c_1}{2-\alpha}\int_{\ep}^t (-\ln
s)^{\frac{\alpha-1}{2-\alpha}}ds.
\endeq
Note that
\neweq{HH0}
\di \int_{0}^t (-\ln
s)^{\frac{\alpha-1}{2-\alpha}}ds<+\infty\qquad\mbox{ and }\;\;
\lim_{t\ri 0^+}\frac{\int_{0}^t (-\ln
s)^{\frac{\alpha-1}{2-\alpha}}ds}{t(-\ln
t)^{\frac{1}{2-\alpha}}}=0.
\endeq
Thus, taking $\ep\ri 0^+$ in
\eq{HH3} we deduce that there exist $c_2>0$ and
$\delta_2\in(0,\delta_1)$ such that
\neweq{HH4}
\di H(t)\leq c_2t(-\ln t)^{\frac{1}{2-\alpha}},\quad\mbox{for all
}\;0<t\leq \delta_2.
\endeq
From \eq{dunus} and \eq{HH4} we obtain
$$\di -H''(t)\geq c_2^{\alpha-1}t^{-1}(-\ln
t)^{\frac{\alpha-1}{2-\alpha}},\quad\mbox{for all
}\;0<t\leq \delta_2.$$ Integrating over $[t,\delta_2]$ in the
above inequality we get
$$\di H'(t)\geq (2-\alpha)c_2^{\alpha-1}
\left[(-\ln t)^{\frac{1}{2-\alpha}}-(-\ln \delta_2)^{\frac{1}{2-\alpha}}\right],\quad\mbox{for all
}\;0<t\leq \delta_2.$$
Therefore, there exist $c_3>0$ and $\delta_3\in(0,\delta_2)$ such that
$$\di H'(t)\geq c_3 (-\ln t)^{\frac{1}{2-\alpha}},
\quad\mbox{for all
}\;0<t\leq \delta_3.$$
With the same arguments as in \eq{HH2}-\eq{HH4}
we obtain $c_4>0$ and
$\delta_4\in(0,\delta_3)$ such that
\neweq{HH5}
\di H(t)\geq c_4t(-\ln t)^{\frac{1}{2-\alpha}},\quad\mbox{for all
}\;0<t\leq \delta_4.
\endeq
The conclusion of (ii) in Theorem \ref{th5} follows now from
\eq{HH4} and \eq{HH5}.

(ii3) Using the fact that $H'(0+)\in(0,+\infty]$ and the
inequality
\eq{HH}, we get the existence of $c>0$ such that
$$\di H(t)>ct,\quad\mbox{for all
}\;0<t\leq b.$$ This yields
$$-H''(t)\leq c^{-\beta} t^{-(\alpha+\beta)},\quad\mbox{for all
}\;0<t\leq b.$$ Since $\alpha+\beta<1,$ it follows that
$H'(0+)<+\infty,$ that is, $H\in C^1[0,b].$ Thus, there exists
$c_1,c_2>0$ such that
\neweq{HH6}
c_1t\leq H(t)\leq c_2t,\quad\mbox{for all }\;0<t\leq b.
\endeq
The conclusion in Theorem \ref{th5} (iii) follows directly from
\eq{HH6} and \eq{asimp1}.

This completes the proof of Theorem \ref{th5}.
\endproof

\noindent {\bf Acknowledgments.} The authors are partially
supported by Programme EGIDE-Brancusi PAI 08915PG between
Université de Picardie Jules Verne and University of
Craiova. This work has been completed while V.R. was visiting the
Université de Picardie Jules Verne in February 2006.
He thanks Prof. O.~Goubet for invitation and for many constructive
discussions.

\end{document}